\documentstyle[12pt,leqno,amssymb,graphics]{amsart}

\newtheorem{thm}{Theorem}[section]
\newtheorem{lemma}[thm]{Lemma}

\theoremstyle{definition}
\newtheorem{dfn}[thm]{Definition}
\theoremstyle{remark}

\begin{document}

\newcommand{\ct}{\cite}
\newcommand{\pr}{\protect\ref}
\newcommand{\su}{\subseteq}
\newcommand{\pa}{{\partial}}
\newcommand{\im}{{Imm(F,\E)}}
\newcommand{\lm}{{\lambda}}
\newcommand{\tc}{{\mathfrak{T}}}
\newcommand{\Q}{{\Bbb Q}}
\newcommand{\hf}{{1 \over 2}}
\newcommand{\CC}{{\mathcal C}}

\newcommand{\R}{{\Bbb R}}
\newcommand{\Z}{{\Bbb Z}}
\newcommand{\E}{{{\Bbb R}^3}}
\newcommand{\C}{{{\Bbb Z}/2}}
\newcommand{\B}{{\Bbb B}}

\newcommand{\FF}{{\mathcal F}}
\newcommand{\A}{{\mathcal{A}}}
\newcommand{\Ai}{{{\mathrm{Arf}}(g^i)}}
\newcommand{\Ar}{{\mathrm{Arf}}}
\newcommand{\r}{{\mathrm{rank}}}
\newcommand{\I}{{\mathrm{Id}}}
\newcommand{\spp}{{\mathrm{supp}}}

\newcommand{\F}{{\mathbf{F}}}
\newcommand{\G}{{\Bbb G}}
\newcommand{\cd}{{\mathrm{codim}}}
\newcommand{\p}{{\psi}}
\newcommand{\hp}{{\Psi}}
\newcommand{\hhp}{\widehat{\Psi}}
\newcommand{\pp}{{\widehat{\psi}}}

\newcommand{\4}{{\mathcal{H}}}
\newcommand{\N}{{\mathcal{N}}}
\newcommand{\ng}{{{\mathcal{N}}_g}}

\newcommand{\hM}{{\widehat{\mathcal{M}}}}
\newcommand{\m}{{\mathcal{M}}_g}
\newcommand{\hm}{{\hM_g}}
\newcommand{\hmi}{{\hM_{g^i}}}
\newcommand{\mi}{{\M_{g^i}}}

\newcommand{\f}{{\overline{f}}}

\newcommand{\hc}{{H_1(F,\C)}}
\newcommand{\ohg}{{O(\hc,g)}}

\newcommand{\ohi}{{O(\hc,g^i)}}
\newcommand{\Tc}{{{\mathcal{T}}_c}}
\newcommand{\SP}{{{\mathcal{S}}_P}}
\newcommand{\Y}{{\mathcal{Y}}}
\newcommand{\T}{{\mathcal{T}}}
\newcommand{\U}{{\mathcal{U}}}
\newcommand{\tF}{{\widetilde{F}}}

\newcommand{\one}{{\widetilde{F_1}}}
\newcommand{\two}{{\widetilde{F_2}}}
\newcommand{\tfk}{{\widetilde{F_k}}}
\newcommand{\th}{{\widetilde{h}}}
\newcommand{\ep}{{\epsilon}}
\newcommand{\wb}{{W^\bot}}
\newcommand{\w}{{\bigcup_i w_i}}
\newcommand{\eep}{{\widehat{\epsilon}}}

\newcounter{numb}

\title{Framings and Projective Framings for 3-Manifolds}
\author{Tahl Nowik}
\address{Department of Mathematics, Bar-Ilan University, 
Ramat-Gan 52900, Israel.}
\email{tahl@@math.biu.ac.il}
\date{December 15, 2003}
\thanks{Partially supported by the Minerva Foundation}

\begin{abstract}
We give an elementary proof of the fact that any orientable 3-manifold admits 
a framing (i.e. is parallelizable) and any non-orientable  3-manifold admits a projective framing.
The proof uses only basic facts about immersions of surfaces in 3-space. 
\end{abstract}

\maketitle

\section{Introduction}

A framing for a smooth $n$-manifold $M$ is a smooth choice of ordered basis $(v_1,\dots,v_n)$ for the
tangent space at each point of $M$. A \emph{projective} framing is a smooth choice of
a \emph{pair} of bases of the form $\{(v_1,\dots,v_n) , (-v_1,\dots,-v_n)\}$.
It has long been established that any orientable  3-manifold admits 
a framing (\ct{s}) and any non-orientable  3-manifold admits a projective framing (\ct{hh}).
The original proofs 
rely on the notion of characteristic classes. 
We present a proof for compact $M$,
which will only use basic facts about immersions of surfaces in 3-space. 

We now present the facts on immersions that we will need.
Denote $\4=(\hf\Z)/(2\Z)$, which is a cyclic group of order 4. 
Let $U$ denote an annulus or Mobius band. There are two regular homotopy classes 
of immersions of $U$ into $\E$. To each such class we attach a 
value in $\4$ as follows: For $U$ an annulus, the regular homotopy class of immersions which includes
an embedding whose image is $S^1 \times [0,1] \su \R^2 \times \R = \E$, will have value $0\in\4$. 
The other class, (containing an embedding differing from the previous embedding by one full twist) will have value 
$1\in\4$.
For $U$ a Mobius band we attach once and for all the value $\hf\in\4$ to one of the classes and the value 
$-\hf\in\4$ 
to the other (which again differ from the previous by one full twist).
Now let $F$ be a closed surface and $i:F \to \E$ an immersion.
We define a map $g^i:\hc \to\4$ as follows:
Given $x \in \hc$ let $c \su F$ be an embedded circle which represents $x$.
Let $U$ be a thin neighborhood of $c$ in $F$, then $U$ is an annulus or Mobius band.
We define $g^i(x)$ to be the value in $\4$ attached above to the immersion $i|_U$.
It has been shown in \ct{p} that $g^i$ is indeed well defined on $\hc$ and satisfies the following property:
\begin{equation}\label{fp}
g^i(x+y) = g^i(x) + g^i(y) + x \cdot y
\end{equation}
where $x \cdot y$ denotes the intersection form on $\hc$. Note $x \cdot y \in \C \su \4$.

We remark that the notation in \ct{p} differs from ours in
that $\4$ is taken there to be $\Z / 4\Z$ rather than $(\hf\Z) / 2\Z$. And so
the numerical value of $g^i$ appearing there is twice the value here, and
our property (\pr{fp}) appearing above is replaced there by 
$g^i(x+y)=g^i(x)+g^i(y)+ 2(x \cdot y)$.

We point out that the above facts about immersions are indeed basic,
in the sense that they may be obtained without knowledge of the Smale-Hirsch Theorem.
E.g. for the classification of regular homotopy classes of immersions of $U$ into $\E$ where $U$ is an 
annulus or Mobius band, one needs to show that the easily constructed map into $\pi_1(SO_3)$ 
is bijective. In general one would use the Smale-Hirsch Theorem, but in this case surjectivity follows
by direct construction of the two immersions, and injectivity follows by familiarity
with the ``belt trick''.

\section{Proof of Theorem}

We will prove the statement for $M$ a closed 3-manifold. The statement for $M$ with boundary
will follow by restricting a framing from the double of $M$. 
So let $M$ be a closed 3-manifold, and let $F \su M$ be a 
Heegard surface, i.e. $F$ splits $M$ into two handlebodies $A,B$ (orientable or non-orientable)
with common boundary $F$. 
Let $a_1,\dots,a_n \su F$ be a system of disjoint circles  
which may be compressed in $A$, reducing $A$ to a ball, and let $b_1,\dots,b_n \su F$ be such system
compressible in $B$. Then
a thin neighborhood $U$ of each of the $a_k$s and $b_k$s is an annulus.
If $M$ is orientable then embed $A$ into $\E$ as in Figure \pr{f1}a.
This induces an embedding $i:F \to \E$ as its boundary. 
If $M$ is non-orientable, then take an immersion $i:F\to\E$ into $\E$ as in Figure \pr{f1}b. 
In both cases we have $g^i(a_k)=0$ for all $1 \leq k \leq n$. (We allow $a_k$ to denote
both the circle in $F$ and its homology class in $\hc$.)

\begin{figure}[t]
\scalebox{0.7}{\includegraphics{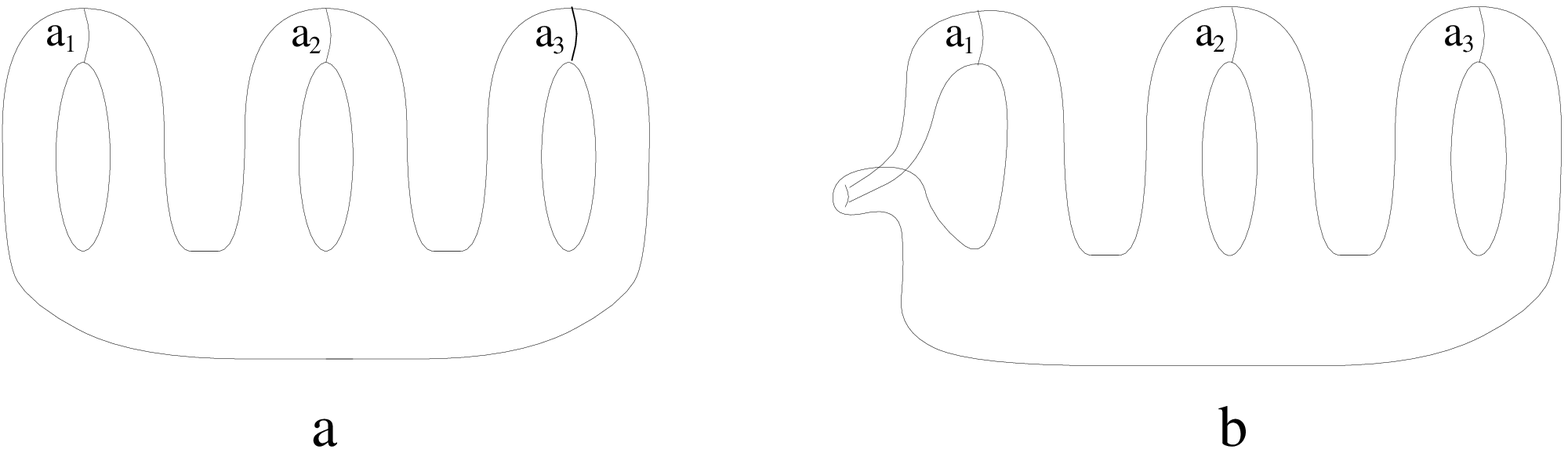}}
\caption{}\label{f1}
\end{figure}

\begin{lemma}\label{dhn}
There exists an $h:F \to F$ which is a composition of Dehn twists along some of the 
$a_k$s, such that $g^i(h(b_k))= 0$ for all $1 \leq k \leq n$.
\end{lemma}

\begin{pf}
For any diffeomorphism $h:F \to F$, define $\phi_h:\hc\to\C$ by
$\phi_h(x) = g^{i\circ h}(x) - g^i(x)$. Note that indeed $\phi_h(x) \in \C$ and by (\pr{fp}) 
$\phi_h$ is a linear functional.
Let $h_k:F \to F$ denote the Dehn twist along $a_k$, then 
the induced homomorphism ${h_k}_* : \hc \to \hc$ is given by 
${h_k}_*(x)=x + (x \cdot a_k)a_k$. We get (by (\pr{fp})) for all $x \in \hc$:
\begin{equation}\label{e2}
\phi_{h_k}(x)= g^{i \circ h_k}(x) - g^i(x) = g^i({h_k}_*(x)) - g^i(x) =
(x \cdot a_k)^2 = x \cdot a_k.
\end{equation}
Let $V_A, V_B \su \hc$ be the subspaces spanned by the $a_k$s and 
$b_k$s respectively, then the intersection form vanishes on $V_A$ and $V_B$. 
It follows by (\pr{fp}) that $g^i|_{V_A}$ is identically 0 and $g^i|_{V_B}$ is a $\C$ valued linear functional. 
By (\pr{e2})
$V_A$ is contained in the kernel of each $\phi_{h_k}$ and since the 
intersection form is non-degenerate on $\hc$, $\phi_{h_1},\dots,\phi_{h_n}$
span all linear functionals on $\hc$ which vanish on $V_A$.
In particular, they span a linear functional $\phi$ vanishing on $V_A$ and satisfying $\phi|_{V_B} = g^i|_{V_B}$.
Now let $h$ be the composition of Dehn twists along those $a_k$s participating
in the presentation of $\phi$ as linear combination of $\phi_{h_1},\dots,\phi_{h_n}$,
then $\phi_h=\phi$ and so for all $1 \leq k \leq n$: 
$g^i(h(b_k))=g^{i \circ h}(b_k)=g^i(b_k)+\phi_h(b_k)=2g^i(b_k)=0$.
\end{pf}

We now change the gluing between $A$ and $B$ using the 
$h:F \to F$ given by Lemma \pr{dhn}. 
This does not change the diffeomorphism type of $M$ since Dehn twists along the $a_k$s may be extended to $A$.
After this change of gluing we have that $g^i(b_k)=0$ for all $1 \leq k \leq n$, and still
$g^i(a_k)=0$ for all $1 \leq k \leq n$.

We will first construct a (projective) framing on $TM|_F$. Let $TF_p \su TM_p$ denote the
tangent spaces of $F$ and $M$ at $p \in F$. Let $n_p \in TM_p$ be the unit vector normal to $TF_p$ and
pointing away from $A$. (Assume we have chosen a Riemannian metric on $M$.)
If $M$ is orientable then there is well defined a unit normal $u_p$ at $i(F) \su \E$ pointing away from the embedding
of $A$. And so $i_*$ together with the choice of normals $n_p$ and $u_p$ induce an isomorphism of $TM_p$ with $T\E_{i(p)}$, with
which we pull back the standard framing of $T\E_{i(p)}$ to $TM_p$.
If $M$ is non-orientable then a normal to $i(F)$ may not be well chosen. let $u_p$ be one of the two unit normals.
We define the isomorphism from $TM_p$ to $T\E_{i(p)}$ using this $u_p$, and pull back the frame in $T\E_{i(p)}$ which 
is obtained from the standard frame by rotating it by $\pi \over 2$ around $u_p$ in the positive sense defined by 
the direction of $u_p$.
When choosing the opposite $u_p$, then the isomorphism will differ by the reflection defined by $u_p$, 
and the frame we pull back from $\E$
will differ by a $\pi$ rotation around $u_p$, and so altogether 
the frame induced in $TM_p$ will differ by 
$-\I$. So such pull back defines a
projective framing on $TM|_F$.

The next step is to extend the framing to a thin neighborhood of the compressing discs of the $a_k$s and $b_k$s
in $A,B$. Let $c$ be one of the circles $a_k$ or $b_k$ and let $U$ be a thin neighborhood of $c$ in $F$, 
then $U$ is an annulus, and since $g^i(c)=0$, there is a regular homotopy $H_t : U \to \E$, $0 \leq t \leq 1$, with
$H_0 = i|_U$ and $H_1$ an embedding with image 
$S^1 \times I \su \R^2 \times \R = \E$ where $I$ denotes $[0,1]$. 

Let $E$ be a neighborhood of the 
compressing disc of $c$ in $A$ or $B$,
so $E = D^2 \times I$ where $U$ corresponds to $\pa D^2 \times I$.
For $0<r\leq 1$ let $E_r , A_r \su E$ denote $D_r \times I, S_r \times I$ where 
$D_r,S_r$ are respectively the disc and circle of radius $r$ in $D^2$. 
For any
$p \in A_r$ take $n_p$ to be the unit normal to $A_r$ pointing inward in $E$.
Now if $M$ is orientable then we follow the increasing time $0 \leq t \leq 1$ of our regular homotopy $H_t$ 
with decreasing radius from $1$
to $\hf$ in $D^2$, and pull back the standard framing of $\E$ to any point $p \in A_r, \hf \leq r \leq 1$
as before, using the normal $n_p$, and normal $u_p$ to the image of the immersion in $\E$, where the side
of the normal is chosen at time 0 to be pointing outward from the embedding of $A$. 
This extends the framing we have constructed in $F$.
At $r=\hf$ we have an embedding of $A_\hf$ as $S^1 \times I \su \R^2 \times \R$, and by additional regular homotopy 
if necessary we may assume the normals $u_p$ are pointing inward. 
This extends to an embedding of $E_\hf$
onto $D^2 \times I \su \R^2 \times \R$. Now if the metric on $M$ is chosen so that this embedding of $E_\hf$ 
into $\E$
is isometric, then the pull back of the standard framing of $\E$ to $A_\hf$ that we have, coincides with the one 
induced by this embedding of $E_\hf$, and so the framing 
can be extended to the whole of $E_{\hf}$.
This completes the construction of the framing in $E$.
We now only need to extend the framing to the remaining 3-balls in $A$ and $B$, but this may
always be done (e.g. since $\pi_2(SO_3)=0$). This completes the proof for the case $M$ is 
orientable.

For non-orientable $M$ we argue as follows: 
Since $U$ is an annulus, it is orientable and so we may choose a unit normal
$u_p$ in $\E$ along $i|_U:U\to\E$. From the projective framing we have defined along $F$ select 
the proper (i.e. non-projective)
framing determined above by this choice of $u_p$. 
We will extend this to a proper framing on $E$ as in the
orientable case, only that now the
framing in $E$ will be defined in three layers, ${2 \over 3} \leq r \leq 1$, ${1\over 3} \leq r \leq {2 \over 3}$,
$0 \leq r \leq {1 \over 3}$. In the outer layer we will realize a rotation of the framing around $u_p$ in the
negative sense, until we have the actual pull back of the framing in $\E$. Then the middle and inner layers
will be defined as in the orientable case, using a regular homotopy of $U$ to a standard annulus in $\E$.
Again we may complete the framing to the remaining balls in $A,B$ since once we have a projective
framing on the boundary of such ball, we may select a proper framing from it (since the ball is orientable),
and then extend the framing to the ball as before.

\end{document}